\documentclass{article}
\usepackage{amssymb}
\usepackage{amsmath}

\usepackage{cite}

\newtheorem{theorem}{Theorem}

\newtheorem{corollary}[theorem]{Corollary}

\newtheorem{definition}[theorem]{Definition}

\newtheorem{lemma}[theorem]{Lemma}

\newtheorem{proposition}[theorem]{Proposition}

\begin{document}

\title{Affine structures on nilpotent contact Lie algebras}
\author{Elisabeth Remm \\
Universit\'{e} de Haute Alsace\\
MULHOUSE 68093}
\date{}
\maketitle

\section{ Introduction} 

\begin{definition}
An affine structure on a $n$-dimensional Lie algebra $\frak{g}$ over $\Bbb{R}$,
is a bilinear map 
\[
\nabla :\frak{g}\times \mathcal{\frak{g}}\rightarrow \frak{g}
\]
satisfying 
\[
\left\{ 
\begin{array}{l}
1)\text{\quad }\nabla \left( X,Y\right) -\nabla \left( Y,X\right) =\left[
X,Y\right]  \\ 
2)\text{\quad }\nabla \left( X,\nabla \left( Y,Z\right) \right) -\nabla
\left( Y,\nabla \left( X,Z\right) \right) =\nabla \left( \left[ X,Y\right]
,Z\right) \quad 
\end{array}
\right. 
\]
for all $X,Y,Z\in \frak{g.}$
\end{definition}

To define an affine structure on a Lie algebra $\frak{g}$ corresponds to give a left invariant affine
flat and torsionfree connection  on a connected Lie group $G$ of associated Lie algebra $\frak{g}.$

\smallskip 

Any Lie algebra equipped with a symplectic form can be equipped with an affine structure. On the other hand there 
exist $(2p+1)$-dimensional Lie algebras  with contact form and no affine structure. But each nilpotent
contact Lie algebra is a one-dimensional central extension of a symplectic algebra. 
The aim of this work is to study how we can extend, under certain conditions, the symplectic stucture.  

\section{Nilpotent Lie algebras with a contact form}

\begin{definition}
Let $\frak{g}$ be an $(2p+1)$-dimensional algebra. A contact form on $\frak{g}$ is a linear form $\omega \neq 0$ of $\frak{g}^*$
such that $\omega \wedge(dw)^p \neq 0.$ In this case $(\frak{g},\omega)$ or $\frak{g}$ is called a contact Lie algebra.
\end{definition} 

\begin{proposition}
\bigskip {\rm ([G])} Let $\frak{g}$ be a contact nilpotent Lie algebra.
Then the center $Z(\frak{g)}$ is one-dimensional.
\end{proposition}

\noindent {\it Proof.} If $\frak{g}$ is $(2p+1)$-dimensional and equipped with a contact form $\omega$,
$\dim Z(\frak{g})\leq 1.$ This follows the fact that if we suppose that $\omega (Z(\frak{g}))=0$ then 
\[
\forall X\in Z(\frak{g})\quad d\omega (X,Y)=-\omega \left[ X,Y\right] =0. 
\]
Thus there exists $X$ such that $\omega (X)=0$ and $X\lrcorner d\omega =0,$ where $\lrcorner$ denotes the inner product. The vector $X$
belongs to the characteristic subspace and $\omega \wedge d\omega
^{p}=0.$ Thus $\omega (Z(\frak{g}))\neq 0$ which proves that $\dim Z(\frak{g}%
)\leq 1.$ If moreover the Lie algebra $\frak{g}$ is nilpotent then 
$\dim Z(\frak{g})=1$ as the center of a nilpotent Lie algebra is never zero. {\rule{.1in}{.1in}}

\begin{corollary}
Let $\frak{g}$ be a contact nilpotent Lie algebra. Then  $\frak{g/}Z\frak{(g)}$ is a symplectic Lie algebra.
\end{corollary}

Thus any contact Lie algebra is a one-dimensional central extension of a symplectic Lie algebra:
\[
0\rightarrow V\rightarrow \frak{g}_{2p+1}\rightarrow (\frak{g}_{2p},\theta
)\rightarrow 0. 
\]

As any symplectic nilpotent Lie algebra can be equipped with an affine structure, we have

\begin{proposition}
Any nilpotent contact Lie algebra is a one-dimensional central extension
of an affine Lie algebra.
\end{proposition}

\section{Affines structures on nilpotent contact Lie algebras}

\subsection{Extension of symplectic affine structure}
Let $(\frak{g}_{2p},\theta )$ be a symplectic nilpotent Lie algebra and $\nabla$ the affine structure coming from this symplectic forme 
that is
\[
\nabla _{X}Y=f(X)Y 
\]
where $f(X)$ is the following endomorphism:

\[
\forall X,Y,Z\in \frak{g\qquad }\theta \left( f(X)(Y),Z\right) =-\theta (Y,%
\left[ X,Z\right] ). 
\]

Let $\widetilde{\frak{g}}$ be the contact Lie algebra defined by the one-dimensional extension 
\[
0\rightarrow V\rightarrow \frak{g}\stackrel{\pi }{\rightarrow } \widetilde{\frak{g}} \rightarrow 0. 
\]
The Lie algebra $\widetilde{\frak{g}}$ identified with $\frak{g\oplus }V$ has the following brackets 
\[
\left[ \left( X,\alpha \right) ,\left( Y,\lambda \right) \right] \widetilde{%
_{\frak{g} }}=\left( \left[ X,Y\right] _{\frak{g} },\theta \left( X,Y\right) \right) 
\]
Let $\widetilde{\nabla }:\widetilde{\frak{g}}\otimes \widetilde{\frak{g}}\rightarrow \widetilde{\frak{g}}$
be an operator satisfying 
\[ 
(*)
\left\{ 
\begin{array}{l}
\widetilde{\nabla }\left( \left( X,0\right) ,\left( Y,0\right) \right)
=\left( \nabla \left( X,Y\right) ,\varphi \left( X,Y\right) \right) \\ 
\widetilde{\nabla }\left( \left( X,0\right) ,\left( 0,\lambda \right)
\right) =\widetilde{\nabla }\left( \left( 0,\lambda \right) ,\left(
X,0\right) \right)
\end{array}
\right. 
\]
where $\varphi $ is a bilinear map on $\frak{g}$
such as
\[
\varphi \left( X,Y\right) -\varphi \left( Y,X\right) =\theta \left(
X,Y\right) . 
\]

\begin{lemma}
the operator $\widetilde{\nabla}$ satisfies the following identity:
$$\widetilde{\nabla }\left( \left( X,\alpha \right) ,\left( Y,\lambda \right)
\right) -\widetilde{\nabla }\left( \left( Y,\lambda \right) ,\left( X,\alpha
\right) \right)=\left[ \left( X,\alpha \right) ,\left( Y,\lambda \right) \right] \widetilde{_{\frak{g} }}$$
\end{lemma}

\noindent {\it Proof.} We have for all $X,Y\in \frak{g}$ and $\quad \lambda ,\mu \in \Bbb{K}$%
\[
\Bbb{
\begin{tabular}{l}
$\widetilde{\nabla }\left( \left( X,\alpha \right) ,\left( Y,\lambda \right)
\right) -\widetilde{\nabla }\left( \left( Y,\lambda \right) ,\left( X,\alpha
\right) \right) $ \\ 
=$\left( \nabla \left( X,Y\right) ,\varphi \left( X,Y\right) \right)
+\lambda \widetilde{\nabla }\left( \left( X,0\right) ,\left( 0,1\right)
\right) +\alpha \widetilde{\nabla }\left( \left( 0,1\right) ,\left(
Y,0\right) \right) $ \\ 
$+\alpha \lambda \widetilde{\nabla }\left( \left( 0,1\right) ,\left(
0,1\right) \right) -\left( \nabla \left( Y,X\right) ,\varphi \left(
Y,X\right) \right) -\alpha \widetilde{\nabla }\left( \left( Y,0\right)
,\left( 0,1\right) \right) $ \\ 
$-\lambda \widetilde{\nabla }\left( \left( 0,1\right) ,\left( X,0\right)
\right) -\lambda \alpha \widetilde{\nabla }\left( \left( 0,1\right) ,\left(
0,1\right) \right) $ \\ 
=$\left( \left[ X,Y\right] _{\frak{g} },\theta \left( X,Y\right) \right) $%
\end{tabular}
} 
\]
But
\[
\left[ \left( X,\alpha \right) ,\left( Y,\lambda \right) \right] \widetilde{%
_{\frak{g} }}=\left( \left[ X,Y\right] _{\frak{g} },\theta \left( X,Y\right) \right), 
\]
which implies that
\[
\widetilde{\nabla }\left( \left( X,\alpha \right) ,\left( Y,\lambda \right)
\right) -\widetilde{\nabla }\left( \left( Y,\lambda \right) ,\left( X,\alpha
\right) \right) =\left[ \left( X,\alpha \right) ,\left( Y,\lambda \right) %
\right] \widetilde{_{\frak{g} }}. \qquad \qquad {\rule{.1in}{.1in}} 
\]
The operator $\widetilde{\nabla }$ is associated to a flat torsionfree connection 
on $\widetilde{\frak{g}}$ .

\medskip

We can note that if $\widetilde{\nabla '}$ is another bilinear map on $\widetilde{\frak{g}}$ such that
$\pi^* \widetilde{\nabla '}=\nabla$, the nullity of the torsion of the linear connection associated to $\widetilde{\nabla '}$
implies that $\widetilde{\nabla '}$ satisfies the same conditions $(*)$. This justifies the choice of the conditions $(*)$.

\bigskip

As we want that $\widetilde{\nabla}$ defines an affine structure on $\frak{g}$, 
we introduce, in order to study the curvature of the linear connection associated to $\widetilde{\nabla}$, the following 
application: 
\[
\begin{tabular}{l}
$C($ $\left( X,\alpha \right) ,\left( Y,\lambda \right) ,(Z,\rho ))$ \\ 
$=\widetilde{\nabla }\left( \left( X,\alpha \right) ,\widetilde{\nabla }%
\left( \left( Y,\lambda \right) ,(Z,\rho )\right) \right) -\widetilde{\nabla 
}\left( \left( Y,\lambda \right) ,\widetilde{\nabla }\left( \left( X,\alpha
\right) ,(Z,\rho )\right) \right) $ \\ 
$-\widetilde{\nabla }\left( \left[ \left( X,\alpha \right) ,\left( Y,\lambda
\right) \right] \widetilde{_{\frak{g} }},(Z,\rho )\right) .$%
\end{tabular}
\]
This gives: 
\[
\begin{array}{l}
C(\left( X,\alpha \right) ,\left( Y,\lambda \right) ,(Z,\rho )) \\ 
=\widetilde{\nabla } ( \left( X,\alpha \right) ,\left( \nabla \left(
Y,Z\right) ,\varphi \left( Y,Z\right) \right) +\rho \widetilde{\nabla }
\left( \left( Y,0\right) ,\left( 0,1\right) \right) +\lambda \widetilde{
\nabla }\left( \left( 0,1\right) ,\left( Z,0\right) \right) + 
\\ 
 \lambda \rho 
\widetilde{\nabla }\left( \left( 0,1\right) ,\left( 0,1\right) \right)
) 
-\widetilde{\nabla }( \left( Y,\lambda \right) ,\left( \nabla \left(
X,Z\right) ,\varphi \left( X,Z\right) \right) +\rho \widetilde{\nabla }%
\left( \left( X,0\right) ,\left( 0,1\right) \right) + \\
\alpha \widetilde{
\nabla }\left( \left( 0,1\right) ,\left( Z,0\right) \right) +\alpha \rho 
\widetilde{\nabla }\left( \left( 0,1\right) ,\left( 0,1\right) \right)
) 
-\widetilde{\nabla }( ( \left[ X,Y\right] _{\mu },\theta (
X,Y) ) ,(Z,\rho )) 
\end{array}
\]

\begin{lemma}
The operator $\widetilde{\nabla}$ satisfies:

\noindent 1)
\[
\begin{tabular}{l}
$C(( X,0) ,( Y,0) ,(Z,0)) 
=( 0,\varphi ( X,\nabla ( Y,Z) ) -\varphi (
Y,\nabla ( X,Z) ) -\varphi ( \left[ X,Y\right] _{\mu
},Z) ) $ \\ 
$+\varphi ( Y,Z) \widetilde{\nabla }( ( X,0)
,( 0,1) ) -\varphi ( X,Z) \widetilde{\nabla }
( ( Y,0) ,( 0,1) )  
-\theta ( X,Y) \widetilde{\nabla }( (Z,0),( 0,1)
) $
\end{tabular}
\]

\noindent 2) \[
\begin{tabular}{l}
$C((X,0),(0,1),(Y,0)) 
=\widetilde{\nabla }\left( \left( X,0\right) ,\widetilde{\nabla }\left(
\left( Y,0\right) ,\left( 0,1\right) \right) \right) $ \\
$-\widetilde{\nabla }%
\left( (\nabla (X,Y),0),(0,1)\right) -\varphi \left( X,Y\right) \widetilde{%
\nabla }\left( \left( 0,1\right) ,\left( 0,1\right) \right) $%
\end{tabular}
\]

\noindent 3)  
\[
C(\left( 0,1\right) ,\left( Y,0\right) ,\left( 0,1\right) )=\widetilde{%
\nabla }\left( \left( 0,1\right) ,\widetilde{\nabla }\left( \left(
Y,0\right) ,\left( 0,1\right) \right) \right) -\widetilde{\nabla }\left(
\left( Y,0\right) ,\widetilde{\nabla }\left( \left( 0,1\right) ,\left(
0,1\right) \right) \right) 
\]

\end{lemma}

This follows directly when we develop the expressions.
\bigskip

\begin{proposition}
If   
\[
C((X,0),(0,1),(Y,0))=0 
\]
then $C($ $\left( X,0\right) ,\left( Y,0\right) ,\left(
0,1\right) )=0.$ 
\end{proposition}

In fact
\[
\begin{tabular}{l}
$C(( X,0) ,( Y,0) ,( 0,1) )
=\widetilde{\nabla }( ( X,0) ,\widetilde{\nabla }(
( Y,0) ,( 0,1) ) ) -\widetilde{\nabla }
( ( Y,0) ,\widetilde{\nabla }( ( X,0)
,( 0,1) ) ) $ \\ 
$\qquad  \qquad \qquad \qquad \qquad \quad - \widetilde{\nabla }( ( \left[ X,Y\right] ,\theta (X,Y)
),( 0,1) ) $ \\ 
$=\widetilde{\nabla }( (\nabla (X,Y),0),(0,1)) +\varphi (
X,Y) \widetilde{\nabla }(( 0,1) ,( 0,1)
) 
-\widetilde{\nabla }( (\nabla (Y,X),0),(0,1)) $ \\
$ \quad -\varphi (
Y,X) \widetilde{\nabla }( ( 0,1) ,( 0,1)
) -\widetilde{\nabla }( ( \left[ X,Y\right] ,\theta
(X,Y) ),( 0,1) ) $ \\ 
$=\widetilde{\nabla }( (\left[ X,Y\right] ,0),(0,1)) +\theta (X,Y)
\widetilde{\nabla }( ( 0,1) ,( 0,1) ) 
-\widetilde{\nabla }( \left[ X,Y\right] ,0),( 0,1)) $ \\
$ \quad -\theta (X,Y)\widetilde{\nabla }( 0,1),( 0,1) ) $ \\ 
$=0.$
\end{tabular}
\]

Let us write some necessary conditions for the application $C$ to be equal to zero.
Let $\pi $ be the canonical projection of $\widetilde{\frak{g}}$ on $\frak{g}$, that is : 
\[
\pi (X,\alpha )=X. 
\]

Let us identify $(X,0)$ with $X$ which permits to consider $\frak{g}$
as a vector subspace of $\widetilde{\frak{g}}$. 
Let us denote $V_{X}$ the vector definied by 
\[
V_{X}=\pi (\widetilde{\nabla }\left( \left( X,0\right) ,\left( 0,1\right)
\right) . 
\]
If $C=0$  we have: 
\[
C(\left( X,0\right) ,\left( Y,0\right) ,(Z,0))=0 
\]
which implies that:
\[
\varphi \left( Y,Z\right) V_{X}-\varphi \left( X,Z\right) V_{Y}-\theta
\left( X,Y\right) V_{Z}=0 
\]
for all $X,Y,Z\in \frak{g}$.

\bigskip

\noindent {\bf Remark.} If $\varphi =0$, we have $\theta =0$, the extension is trivial and we came out of the symplectic case.
The operator $\widetilde{\nabla }$ define by 
\begin{eqnarray*}
\widetilde{\nabla }\left( \left( X,0\right) ,\left( Y,0\right) \right)
&=&(\nabla \left( X,Y\right) ,0) \\
\widetilde{\nabla }\left( \left( X,0\right) ,\left( 0,\lambda \right)
\right) &=&0 \\
\widetilde{\nabla }\left( \left( 0,\mu \right) ,\left( 0,\lambda \right)
\right) &=&0
\end{eqnarray*}
give an affine structure on $\widetilde{\frak{g}}=\frak{g} \oplus \mathbb{R}$ which is a direct sum of ideals 

\bigskip

We must then suppose that $\varphi \neq 0.$

\subsection{Case $\protect\varphi =\frac{\protect\theta }{2}.$}

Then we have
\[
\theta \left( \left[ X,Y\right] ,Z\right) =-\theta \left( Y,\nabla
(X,Z)\right) 
\]
and we deduce from the lemma  the following relations: 
\begin{eqnarray*}
&&C(\left( X,0\right) ,\left( Y,0\right) ,(Z,0)) \\
&=&-\frac{1}{2}\left( 0,\theta \left( \left[ X,Y\right] _{\mu },Z\right)
\right) +\frac{1}{2}\theta \left( Y,Z\right) \widetilde{\nabla }\left(
\left( X,0\right) ,\left( 0,1\right) \right)  \\
&&-\frac{1}{2}\theta \left( X,Z\right) \widetilde{\nabla }\left( \left(
Y,0\right) ,\left( 0,1\right) \right) -\theta \left( X,Y\right) \widetilde{%
\nabla }\left( (Z,0),\left( 0,1\right) \right).
\end{eqnarray*}
For all $X \in \frak{g}$ we define $a_X \in \mathbb{R}$ by
$$
\widetilde{\nabla}((X,0),(0,1))=(V_X,a_X).
$$
The nullity of the curvature tensor implies that

\[
(**)
\left\{ 
\begin{array}{l}
\frac{1}{2}\theta \left( Y,Z\right) V_{X}-\frac{1}{2}\theta \left(
X,Z\right) V_{Y}-\theta \left( X,Y\right) V_{Z}=0 \\ 
\theta \left( \left[ X,Y\right] _{\mu },Z\right) +\theta \left( Y,Z\right)
a_{X}-\theta \left( X,Z\right) a_{Y}-2\theta \left( X,Y\right) a_{Z}=0
\end{array}
\right. 
\]
As $\theta $ is of maximal rank, for all $X\in \frak{g}$, there exists 
$Y,Z\in \frak{g}$ such as $\theta \left( X,Z\right) =0=\theta \left(
X,Y\right) $ and $\theta \left( Y,Z\right) =1$. The first of the relations $(**)$
implies that $V_{X}=0.$ Thus 
\[
\widetilde{\nabla }\left( \left( X,0\right) ,\left( 0,1\right) \right)
=(0,a_{X}) 
\]
for all $X\in \frak{g.}$ Then we have 
\[
\widetilde{\nabla }\left( \left( X,0\right) ,\widetilde{\nabla }\left(
\left( Y,0\right) ,\left( 0,1\right) \right) \right) =\widetilde{\nabla }%
\left( \left( X,0\right) ,\left( 0,a_{Y}\right) \right) =a_{X}a_{Y}\left(
0,1\right) 
\]
and

\[
C((X,0),(0,1),(Y,0))=0 
\]
implies that 
\[
a_{X}a_{Y}\left( 0,1\right) =a_{\nabla (X,Y)}\left( 0,1\right) +\frac{1}{2}%
\theta \left( X,Y\right) \widetilde{\nabla }\left( \left( 0,1\right) ,\left(
0,1\right) \right). 
\]
Let us chose $X,Y\in \frak{g}$ such that $\theta \left( X,Y\right) =1.$ We can deduce
\[
\pi (\widetilde{\nabla }\left( \left( 0,1\right) ,\left( 0,1\right) \right)
)=0\frak{\ } 
\]
and then 
\[
\widetilde{\nabla }\left( \left( 0,1\right) ,\left( 0,1\right) \right)
=(0,a), \qquad a \in \mathbb{R}. 
\]
Moreover 
\begin{eqnarray*}
a_{X}a_{Y} &=&a_{\nabla (X,Y)}+\frac{1}{2}\theta \left( X,Y\right) a \\
a_{X}a_{Y} &=&a_{\nabla (Y,X)}-\frac{1}{2}\theta \left( X,Y\right) a
\end{eqnarray*}
which gives 
\begin{eqnarray*}
\theta \left( X,Y\right) a &=&-a_{\nabla (X,Y)}+a_{\nabla (X,Y)} \\
&=&-a_{[X,Y].}
\end{eqnarray*}
Let us denote the linear form $\alpha $ de $\frak{g}^{\ast }$ defined by
$\alpha (X)=a_{X}.$ If $a\neq 0$, then $\theta \left( X,Y\right) =\frac{1}{a}%
d\alpha (X,Y).$ The symplectic cocycle $\theta $ is then exact. But on
any nilpotent Lie algebra, the class of linear forms is
odd ([G]). We deduce that it can not exist exact symplectic form on $\frak{g}$ 
and the previous equality can not be true. If $a=0,$ 
\[
\widetilde{\nabla }\left( \left( 0,1\right) ,\left( 0,1\right) \right) =0 
\]
and 
\[
a_{[X,Y]}=0. 
\]
Then the application 
\[
\alpha :\frak{g\rightarrow }\Bbb{R} 
\]
given by $\alpha (X)=a_{X}$ defines a one-dimensional linear representation. 

\begin{proposition}
The application $\widetilde{\nabla }$ on the contact Lie algebra $\widetilde{\frak{g}}$ $\frak{\ }$
given by 
\[
\left\{ 
\begin{array}{l}
\widetilde{\nabla }\left( \left( X,0\right) ,\left( Y,0\right) \right)
=\left( \nabla \left( X,Y\right) ,1/2\theta \left( X,Y\right) \right) \\ 
\widetilde{\nabla }\left( \left( X,0\right) ,\left( 0,\lambda \right)
\right) =\widetilde{\nabla }\left( \left( 0,\lambda \right) ,\left(
X,0\right) \right)
\end{array}
\right. 
\]
where $\nabla $ is the affine strucutre on $\frak{g}=\frac{\widetilde{\frak{g}}}{Z(\widetilde{\frak{g}})}$ 
associated to the symplectic structure $\theta $ is an affine structure if there exits a one-dimensional representation
of $\frak{g}$ such that
\[
\theta \left( \left[ X,Y\right] _{\mu },Z\right) +\theta \left( Y,Z\right)
a_{X}-\theta \left( X,Z\right) a_{Y}-2\theta \left( X,Y\right) a_{Z}=0 
\]
for all $X,Y,Z\in \frak{g.}$
\end{proposition}

\subsection{General case }

We still suppose  
\[
\theta \left( \left[ X,Y\right] ,Z\right) =-\theta \left( Y,\nabla
(X,Z)\right) 
\]
We saw that the nullity of the curvature of the connection associated to the bilinear application $\widetilde{\nabla}$ implies:

\[
\varphi \left( Y,Z\right) V_{X}-\varphi \left( X,Z\right) V_{Y}-\theta
\left( X,Y\right) V_{Z}=0 
\]
for all $X,Y,Z\in \frak{g}$.

Let us suppose that we have $X\in \frak{g}$ such that $V_{X}\neq 0.$ Let us take a vector $X$ satisfying this property.
The orthogonal space for $\theta $ of the space $\Bbb{R}\left\{ X\right\} $ is of codimension $1$. 
Let $Y\in \Bbb{R}\left\{ X\right\} ^{\bot }$ . Then $\theta \left( X,Y\right) =0$ and

\[
\varphi \left( Y,Z\right) V_{X}=\varphi \left( X,Z\right) V_{Y}, \quad Z\in 
\frak{g} 
\]
If we can find a vector $Z\in $ $\frak{g}$ such as $\varphi \left(Y,Z\right) \neq 0$,
we have that $V_{Y}\neq 0$ and in this case the non zero vectors $V_{X}$ and $V_{Y}$
are colinear. Let us assume that 
\[
V_{X}=\lambda _{X,Y}V_{Y} 
\]
then
\[
\lambda _{X,Y}\varphi \left( Y,Z\right) =\varphi \left( X,Z\right) 
\]
for all $Z\in \frak{g}$. Let us consider a basis $\left(
X_{1},...,X_{2m}\right)$ of $\frak{g}$ in which the matrix of $\theta $ is reduced to the following form: 
\[
\left( 
\begin{array}{ccc}
\begin{array}{cc}
0 & 1 \\ 
-1 & 0
\end{array}
&  & 0 \\ 
& \ddots &  \\ 
0 &  & 
\begin{array}{cc}
0 & 1 \\ 
-1 & 0
\end{array}
\end{array}
\right) 
\]
This shows that the matrix of $\varphi $ is at least of rank 2 because $\varphi
\left( X,Y\right) -\varphi \left( Y,X\right) =\theta \left( X,Y\right).$
For instance, if $\frak{g}$ is 4-dimensional, the matrix of $\varphi $ would be of the following form:

\[
\left( 
\begin{array}{cccc}
\alpha _{1} & \alpha _{2} & \alpha _{3} & \alpha _{4} \\ 
\alpha _{2}-1 & \alpha _{6} & \alpha _{7} & \alpha _{8} \\ 
\alpha _{3} & \alpha _{7} & \alpha _{11} & \alpha _{12} \\ 
\alpha _{4} & \alpha _{8} & \alpha _{12}-1 & \alpha _{16}
\end{array}
\right) 
\]
We can suppose that $\varphi \left( X_{1},\cdot \right) \neq 0$ (in other case, there would exist a vector
$X_{i}$ such that $\varphi \left( X_{i},Z\right) \neq 0$
since $\varphi $ is of non zero rank and it is sufficient to take a basis
adapted to $\theta $ with $X_{i}$ as first vector). 
Let $Y$ belong to vector space  $\left\{ X_{3},...,X_{2n}\right\} $; 
we have that $\theta \left(X_{1},Y\right) =0$
\[
\varphi \left( X_{1},Z\right) =\lambda _{X_{1},Y}\varphi \left( Y,Z\right)
\qquad \mbox{for all} \ Z\in \frak{g} 
\]
and
\[
\lambda _{X_{1},Y}\neq 0 
\]
Therefore the $(2m-2)$ last columns of the matrix of $\varphi $ (associated to the chosen basis)
are proportional to the first one. But we also have $\theta \left( X_{2},X_{3}\right) =0$ and

\[
\varphi \left( X_{2},Z\right) =\lambda _{X_{2},X_{3}}\varphi \left(
X_{3},Z\right) 
\]
Similary we show that $\theta \left( X_{1},X_{3}\right) =0$ and 
\[
\varphi \left( X_{1},Z\right) =\lambda _{X_{1},X_{3}}\varphi \left(
X_{3},Z\right) \qquad \mbox{for all} \ Z\in \frak{g}. 
\]
We deduce that $\varphi \left( X_{2},Z\right) =\lambda \varphi \left( X_{1},Z\right),$
which implies that $\varphi $ is of rank 1. This is impossible and then we have $V_{X}=0.$

\begin{proposition}
Let $\widetilde{\frak{g}}$ be a contact nilpotent Lie algebra. If the affine structure $\nabla $ which is 
defined by a symplectic cocycle on $\frak{\bar{g}=}\frac{\widetilde{\frak{g}}}{Z(\widetilde{\frak{g}})}$
can be extended to an affine structure $\widetilde{\nabla }$ on $\widetilde{\frak{g}}$, we have:  
\[
\pi (\widetilde{\nabla }\left( X,T\right) )=0 
\]
for all X$\in \frak{g}$ and $T\in Z(\widetilde{\frak{g}}).$
\end{proposition}

We have that for all vector $X$ in $\frak{g}$, $V_{X}=0$ and  
$\widetilde{\nabla }\left(\left( X,0\right) ,\left( 0,1\right) \right) =\left( 0,a_{X}\right) .$ Then the
equality 
$C(\left( X,0\right) ,\left( Y,0\right) ,(Z,0))=0$ 
implies that
\[
\begin{tabular}{l}
$\varphi \left( X,\nabla \left( Y,Z\right) \right) -\varphi \left( Y,\nabla
\left( X,Z\right) \right) -\varphi \left( \left[ X,Y\right] _{\mu },Z\right) 
$ \\ 
$=-a_{X}\varphi \left( Y,Z\right) +a_{Y}\varphi \left( X,Z\right)
+a_{Z}\theta \left( X,Y\right) $
\end{tabular}
\]
Similarly $C((X,0),(0,1),(Y,0))=0$ implies that
\[
\widetilde{\nabla }\left( \left( X,0\right) ,\left( 0,a_{Y}\right) \right)
-\left( 0,a_{\nabla (X,Y)}\right) -\varphi \left( X,Y\right) \widetilde{%
\nabla }\left( \left( 0,1\right) ,\left( 0,1\right) \right) =0. 
\]
This gives the following equation
\[
\varphi \left( X,Y\right) \widetilde{\nabla }\left( \left( 0,1\right)
,\left( 0,1\right) \right) =(a_{Y}a_{X}-a_{\nabla (X,Y)})(0,1) 
\]
and
\[
\varphi \left( Y,X\right) \widetilde{\nabla }\left( \left( 0,1\right)
,\left( 0,1\right) \right) =(a_{Y}a_{X}-a_{\nabla (Y,X)})(0,1) 
\]
if we permute the vectors $X$ and $Y$.
We combine this two equations to obtain:
\begin{eqnarray*}
\theta (X,Y)\widetilde{\nabla }\left( \left( 0,1\right) ,\left( 0,1\right)
\right) &=&(a_{\nabla (Y,X)}-a_{\nabla (X,Y)})(0,1) \\
&=&a_{[X,Y]}(0,1).
\end{eqnarray*}
This shows in particular that $\widetilde{\nabla }\left( \left( 0,1\right)
,\left( 0,1\right) \right) =\rho (0,1)$ and 
\[
\rho \theta (X,Y)=a_{[X,Y]} .
\]
Finally $C(\left( 0,1\right) ,\left( Y,0\right) ,\left( 0,1\right) )=0$
implies

\[
a_{Y}\widetilde{\nabla }\left( \left( 0,1\right) ,\left( 0,1\right) \right) =%
\widetilde{\nabla }\left( \left( Y,0\right) ,\widetilde{\nabla }\left(
\left( 0,1\right) ,\left( 0,1\right) \right) \right) 
\]
thus 
\[
a_{Y}\widetilde{\nabla }\left( \left( 0,1\right) ,\left( 0,1\right) \right)
=\rho \left( 0,a_{Y}\right). 
\]
This last equation is already satisfied.

Then let us suppose $\rho \neq 0.$ In this case $a_{[X,Y]}\neq 0$ when $\theta (X,Y)\neq 0$. 
Let us take $X$ in $Z(\frak{g}).$ As $\theta $ is of maximal rank, there is one $Y$ such that $\theta (X,Y)\neq 0.$
But $[X,Y]=0$ implies $a_{[X,Y]}=0$. This leads to contradiction.

\medskip

\noindent {\bf Conclusion .} As $\rho=0$ we have that $\widetilde{\nabla }\left( \left(
0,1\right) ,\left( 0,1\right) \right) =0.$ Then $a_{[X,Y]}=0$ and the application 
$\alpha :\frak{g\rightarrow }\Bbb{R}$ defined by $\alpha(X)=a_{X}$ 
gives an one-dimensional linear representation of $\frak{g}$. We deduce

\begin{theorem}
Let $\alpha :\frak{g\rightarrow }\Bbb{R}$ be a one-dimensional linear representation of $\frak{g}$.

When $\alpha $ is the trivial representation, $\widetilde{\nabla }$
is an affine structure if and only if

1) $\widetilde{\nabla }\left( U,\left( 0,1\right) \right) =0$ for all $U\in \widetilde{\frak{g}}$.

2) $\varphi $ satisfies $\varphi \left( X,\nabla \left( Y,Z\right) \right)
-\varphi \left( Y,\nabla \left( X,Z\right) \right) -\varphi \left( \left[ X,Y%
\right] _{\mu },Z\right) =0$, i.e. if it is a 2-cocycle for cohomology of the Vinberg algebra associated to $\nabla $ with values 
in a trivial module.

When $\alpha $ is a non-trivial representation, $\widetilde{\nabla }$ is an affine structure if and only if

1) $\widetilde{\nabla }\left( \left( 0,1\right) ,\left( 0,1\right) \right)
=0 $ , $\widetilde{\nabla }\left( \left( X,0\right) ,\left( 0,1\right)
\right) =0$ for all $X\in Ker \alpha . $

2) $\varphi \left( X,\nabla \left( Y,Z\right) \right) -\varphi \left(
Y,\nabla \left( X,Z\right) \right) -\varphi \left( \left[ X,Y\right] _{\mu
},Z\right) =\alpha (Z)\theta (X,Y)$ for all $X,Y\in Ker$($\alpha ).$
\end{theorem}

\bigskip

\noindent {\bf References}

\bigskip

\noindent [Au\,] Auslander L., {\it The structure of complet locally affine manifolds}. 
Topology {\bf 3} 1964 suppl.1, 131-139. 

\smallskip

\noindent [Be\,] Benoist Y., {\it Une nilvari\'{e}t\'{e} non affine}.
J.Diff.Geom.,   {\bf{41}}, (1995), 21-52.

\smallskip

\noindent [Bu\,] Burde D., {\it Affine structures on nilmanifolds}. 
Int. J. of Math,  {\bf 7} (1996), 599-616.

\smallskip 

\noindent [D-H\,] Dekimpe K., Hartl M., {\it Affine structures on 4-step nilpotent Lie algebras}, 
J. Pure Appl. Algebra {\bf 120} (1997), no.1, 19-37. 

\smallskip

\noindent [F-G\,] Fried D., Goldman W., {\it Three dimensional affine
crystallographic groups.} Adv. Math.,  {\bf \ 47}, (1983), 1-49. 

\smallskip

\noindent [G\,] Goze M., {\it Sur la classe des formes et syst\`{e}mes invariants \`a gauche sur un groupe de Lie.}
CRAS Paris S\'er. A-B {\bf 283} (1976), no.7, Aiii, A499-A502.

\smallskip

\noindent [G-K\,] Goze M., Khakimdjanov Y., {\it Nilpotent Lie algebras}.
Kluwer editor, 1995.

\smallskip

\noindent [G-R\,] Goze M., Remm E., {\it Affine structures on abelian Lie
algebras}, Linear Algebra and its Applications, {\bf 360} (2003), 215-230.

\smallskip

\noindent [H\,] Helmstetter J., {\it Radical d'une alg\`{e}bre sym\'{e}trique
\`{a} gauche.} Ann. Inst. Fourier,  {\bf 29} (1979), 17-35.

\smallskip

\noindent [Ku\,] Kuiper N., {\it Sur les surfaces localement affines}. Colloque
G\'{e}om\'{e}trie diff\'{e}rentielle Strasbourg, (1953), 79-87.

\smallskip

\noindent [Ma\,] Malcev A., {\it Commutative subalgebras of semi-simple Lie algebras}. 
Bull. Acad. Sci. URSS. S\'er. Math. [Izvestia Akad. Nauk SSSR] {\bf 9}, (1945), 291-300.    

\smallskip

\noindent [R\,] Remm E., {\it Structures affines sur les alg\`ebres de Lie et op\'erades Lie-admissibles}, Thesis, december 2001. 

\smallskip

\noindent [S\,] Scheuneman J., {\it affine structures on three-step nilpotent Lie algebras}, 
Proc. Amer. Math. Soc. {\bf 46} (1974), 451-454.

$$
$$
$$
$$
$$
$$

\end{document}